\documentclass[fleqn]{mat01}
\usepackage{times,mathtimy,amssymb,latexsym,amsbsy}
\begin{document}

\setcounter{page}{93}
\firstpage{93}

\newtheorem{theore}{Theorem}
\renewcommand\thetheore{\arabic{section}.\arabic{theore}}
\newtheorem{theor}{\bf Theorem}
\newtheorem{propo}[theore]{\rm PROPOSITION}

\def\definit{\trivlist\item[\hskip\labelsep{{\rm DEFINITION}}]}
\def\claim{\trivlist\item[\hskip\labelsep{{\it Claim.}}]}
\def\profth{\trivlist\item[\hskip\labelsep{{\it Proof of Theorem $1$ for ${S^n}$.}}]}
\def\profthe{\trivlist\item[\hskip\labelsep{{\it Proof of Theorem $2$ for ${S^n}$.}}]}
\def\rempro{\trivlist\item[\hskip\labelsep{{\it Remark on proofs of Theorem $1$ and Theorem $2$ for ${\Hy^n}$.}}]}

\renewcommand{\theequation}{\thesection\arabic{equation}}

\font\german=eufm10 scaled 1200
\font\Bbb=msbm10 scaled 1000

\newcommand{\E}{{\mathbb E}}
\newcommand{\C}{\mathcal {C}}
\newcommand{\N}{{\mathbb N}}
\newcommand{\R}{{\mathbb R}}
\newcommand{\Hy}{{\mathbb H}}
\newcommand{\Om}{\Omega}
\newcommand{\La}{\Delta}
\newcommand{\n}{\nabla}
\newcommand{\Ho}{H^1(\Om)}
\newcommand{\Ht}{H^2(\Om)}
\newcommand{\Hz}{H_0^1(\Om)}
\newcommand{\smb}{\C^\infty(\bar{\Om})}
\newcommand{\smc}{\C_0^\infty(\Om)}
\newcommand{\Lt}{L^2(\Om)}
\newcommand{\lm}{\lambda_1(\Om)}
\newcommand{\F}{\mathcal{F}}

\title{On two functionals connected to the Laplacian in a class of
doubly connected domains in space-forms}

\markboth{M~H~C~Anisa and A~R~Aithal}{Doubly connected domains in
space-forms}

\author{M~H~C~ANISA and A~R~AITHAL}

\address{Department of Mathematics, University of Mumbai,
Mumbai~400~098, India\\
\noindent E-mail: anisa@sankhya.mu.ac.in; aithal@math.mu.ac.in}

\volume{115}

\mon{February}

\parts{1}

\pubyear{2005}

\Date{MS received 7 September 2004; revised 15 December 2004}

\begin{abstract}
Let $B_1$ be a ball of radius $r_1$ in $S^n(\Hy^n)$, and let $B_0$ be
a smaller ball of radius $r_0$ such that $\overline{B_0}\subset B_1$.
For $S^n$ we consider $r_1< \pi$. Let $u$ be a solution of the problem
$-\La u =1$ in $\Om := B_1\setminus \overline{B_0}$ vanishing on the
boundary. It is shown that the associated functional $J(\Om)$ is minimal
if and only if the balls are concentric. It is also shown that the first
Dirichlet eigenvalue of the Laplacian on $\Om$ is maximal if and only if
the balls are concentric.
\end{abstract}


\maketitle

\section{Introduction}

Let $(M, g)$ be a Riemannian manifold and let $D$ denote the
Levi--Civita connection of $(M, g)$. For a smooth vector field $X$ on
$M$ the divergence ${\rm div}(X)$ is defined as ${\rm trace}(DX)$. For a
smooth function $f\hbox{:}\ M \longrightarrow \R$, the gradient $\n f$
is defined by $g(\n f(p), v) = {\rm d}f(p)(v)$ ($p \in M, \; v \in
T_pM$) and the Laplace--Beltrami operator $\La$ is defined by $\La f =
{\rm div}(\n f)$. Further, $\n^2 f$ denotes the Hessian of $f$. Throughout
this paper, $\omega$ and ${\rm d}V$ denote the volume element of $(M,
g)$.

Let $\Om \subset M$ be a domain such that $\bar{\Om}$ is a smooth
compact submanifold of $M$. The Sobolev space $\Ho$ is defined as the
closure of $\smb$ (the space of real valued smooth functions on
$\bar{\Om}$) with respect to the Sobolev norm
\begin{equation*}
\|f\|_{\Ho} = \left( \int_\Om \{ f^2 + \|\n f\|^2 \}\, {\rm d}V
\right)^{1/2} \quad (f \in \smb).
\end{equation*}
The closure of $\smc$ (the space of real valued smooth functions on
$\Om$ having compact support in $\Om$) in $\Ho$ is denoted by $\Hz$. The
Sobolev space $\Ht$ is defined as the closure of $\smb$ with respect to
the Sobolev norm
\begin{equation*}
\|f\|_{\Ht} = \left( \int_\Om \{ f^2 + \|\n f\|^2 + \| \n^2
f\|^2\}\,{\rm d}V\right)^{1/2} \quad (f \in \smb).
\end{equation*}
These spaces are Hilbert spaces with the corresponding norms.

Consider the Dirichlet boundary value problem on $\Om$:
\begin{equation}
\left. \begin{array}{r@{\ }l@{\quad}ll}
-\La u &= 1 &\hbox{on} &\Om,\\[.2pc]
u &= 0 &\hbox{on} &\partial \Om.
\end{array} \right\}
\end{equation}
Let $u \in \Hz$ be the unique weak solution of problem (1.1). By
Theorem~4.8, p.~105 of \cite{1}, $u \in \smb$.

Consider the following eigenvalue problem on $\Om$:
\begin{equation}
\left. \begin{array}{r@{\ }l@{\quad}ll}
-\La u &= \lambda u &\hbox{on} &\Om,\\[.2pc]
u &= 0 &\hbox{on} &\partial \Om.
\end{array} \right\}
\end{equation}
The eigenvalues of the positive Laplace--Beltrami operator $-\La = -{\rm
div}(\n f)$ are strictly positive. The eigenfunctions corresponding to
the first eigenvalue $\lambda_1$ are proportional to each other. They
belong to $\smb$ and they are either strictly positive or strictly
negative on $\Om$. Moreover,
\begin{equation*}
\lambda_1 = \inf \,\{ \, \|\n \phi\|_{\Lt}^2 \,|\; \phi \in \Hz,\,
\|\phi\|_{\Lt}^2 = 1\}
\end{equation*}
(cf. \cite{1}, Theorem~4.4, p.~102). Let $y := y(\Om) \in \smb$ be the
unique solution of problem (1.1). Let $y_1 := y_1(\Om)$ be the unique
solution of problem (1.2), corresponding to the first eigenvalue
$\lambda_1 := \lm$, characterized by
\begin{equation*}
y_1 > 0 \quad \hbox{on }\, \Om \quad \hbox{and} \quad \int_\Om y_1^2\,
{\rm d}V = 1.
\end{equation*}

The aim of this paper is to prove the main results of \cite{3} for
simply connected spherical and hyperbolic space-forms.

Consider the unit sphere $S^n = \{(x_1, x_2, \ldots, x_{n+1}) \in
\R^{n+1}|\sum_{i=1}^{n+1} x_i^2 = 1\}$ with induced Riemannian
metric $\left<\,,\, \right>$ from the Euclidean space $\R^{n+1}$.
Also consider the hyperbolic space $\Hy^n= \{(x_1, x_2, \ldots,
x_{n+1}) \in \R^{n+1}|\sum_{i=1}^{n} x_i^2 -x_{n+1}^2= -1$ and
$x_{n+1}>0\}$ with the Riemannian metric induced from the
quadratic form $\left(x, y\right):=\sum_{i=1}^{n} x_i y_i -
x_{n+1}y_{n+1}$, where $x=(x_1, x_2, \ldots, x_{n+1})$ and
$y=(y_1, y_2, \ldots, y_{n+1})$.

Fix $0<r_0<r_1$. We choose $r_1 < \pi $ for the case of $S^n$. Let
$B_1$ be any ball of radius $r_1$ in $S^n(\Hy^n)$ and $B_0$ be any
ball of radius $r_0$ such that $\overline{B_0}\subset B_1$.
Consider the family $\F = \{ B_1 \setminus \overline{B_0}\}$ of
domains in $S^n(\Hy^n)$. We study the extrema of the following
functionals:
\begin{align*}
J(\Om) &= -\int_\Om \{\, \|\n y(\Om)\|^2 - 2y(\Om) \} \, {\rm d}V,\tag{1}\\[.2pc]
J_1(\Om) &= -\int_\Om \{\, \|\n y_1(\Om)\|^2 - 2\lm [y_1(\Om)]^2 \} \,
{\rm d}V\tag{2}
\end{align*}
on $\F$, associated to problems (1.1) and (1.2) respectively. Note
here that the functionals $J$ and $J_1$ are nothing but negative
of the energy functional $\int_\Om \|\n y(\Om)\|^2\,{\rm d}V$ and
the Dirichlet eigenvalue $\lambda_1$, respectively.

We state our main results: Put $\Om_0 = B(p,\,r_1) \setminus
\overline{B(p, r_0)}$ for any fixed $p \in S^n(\Hy^n)$.

\begin{theor}[\!]
The functional $J(\Om)$ on $\F$ assumes minimum at $\Om$ if and
only if $\Om = \Om_0${\rm ,} i.e.{\rm ,} when the balls are
concentric.
\end{theor}

\begin{theor}[\!]
The functional $J_1(\Om)$ on $\F$ assumes maximum at $\Om$ if and
only if $\Om = \Om_0${\rm ,} i.e.{\rm ,} when the balls are
concentric.
\end{theor}

In \S\S2 and 3, following \cite{5}, we develop the `shape
calculus' for Riemannian manifolds for the stationary problem
(1.1) and the eigenvalue problem (1.2) respectively. In \S4, we
prove Theorems~1 and 2 for $S^n$, and make the necessary remarks
to carry out the proofs of Theorems~1 and 2 for $\Hy^n$.

\section{Shape calculus for the stationary problem}

\setcounter{equation}{0}

Let $V$ be a smooth vector field on $M$ having compact support. Let
$\Phi\hbox{:}\ \R \times M \longrightarrow M$ be the smooth flow for
$V$. For each $t \in \R$, denote $\Phi(t, x)$ by $\Phi_t(x)$ ($x \in
M$). Let $\Om$ be an open subset of $M$ such that $\bar{\Om}$ is a
smooth compact submanifold of $M$. Put $\Om_t := \Phi_t(\Om)$ ($t \in
\R$).

Let $\mathcal{D}$ be a domain in $M$ such that supp$\,V\subset
\mathcal{D}$. Fix $f \in \C^\infty(\mathcal{D})$. Consider the Dirichlet
boundary value problem on $\Om_t$:
\begin{equation}
\left. \begin{array}{r@{\ }l@{\quad}ll}
\La u &= f &\hbox{on} &\Om_t,\\[.2pc]
u &= 0 &\hbox{on} &\partial \Om_t.
\end{array} \right\}
\end{equation}
Let $y_t \in \C^\infty(\bar{\Om}_t)$ be the unique solution of problem
(2.1) (cf. \cite{1}, Theorem~4.8, p.~105). Throughout this section
$y:=y(\Om)$ denotes the unique solution of (2.1) for $t=0$.

Denote $y_t\circ \Phi_t \arrowvert_\Om$ by $y^t$ ($t \in \R$).

\vspace{.4pc}
\begin{propo}$\left.\right.$\vspace{.5pc}

\noindent The map $t \longmapsto y^t$ is a $\C^1$-curve in $\Ht \cap
\Hz$ from a neighbourhood of $0$ in $\R$.\vspace{.4pc}
\end{propo}

\begin{proof}
By problem (2.1), for each $t \in \R$, $y_t$ satisfies the equation
\begin{equation*}
\int_{\Om_t} g(\n y_t,\n \psi)\; {\rm d}V = -\int_{\Om_t} f \psi\; {\rm
d}V \quad \forall \,\psi \in \C_0^\infty(\Om_t).\tag{3}
\end{equation*}
There exists smooth function $\gamma_t\hbox{:}\ M \longrightarrow
(0, \infty)$ such that $\Phi_t^* \omega = \gamma_t \omega$ (here,
$\omega := {\rm d}V$, the volume element of $(M,g)$). Put $B_t :=
(D\Phi_t)^{-1}$, $B_t^* =$ transpose of $B_t$ (i.e., $g(B_t(x)v,
w)= g(v, B_t^*(x)w)\quad \forall v \in T_x\Om_t,\, w \in
T_{x^\prime}\Om$, where $x^\prime:= \Phi_t^{-1}(x)$) and $A_t :=
\gamma_t B_t B_t^*$. By the change of variable $\Phi_t\hbox{:}\
\Om \longrightarrow \Om_t$, eq.~(3) can be re-written as
\begin{equation*}
\int_\Om - {\rm div}(A_t \n (y_t \circ \Phi_t))\; \psi \circ
\Phi_t \; {\rm d}V = -\int_\Om f \circ \Phi_t\psi \circ \Phi_t
\gamma_t\; {\rm d}V.
\end{equation*}
Therefore, $y^t := y_t \circ \Phi_t\hbox{:}\ \Om \longrightarrow \R$
satisfies
\begin{equation}
\left. \begin{array}{r@{\ }l@{\quad}ll}
-{\rm div}(A_t \n y^t) + f \circ \Phi_t \gamma_t  &=  0 &\hbox{on} &\Om,\\[.2pc]
y^t &= 0 &\hbox{on} &\partial \Om.
\end{array} \right\}
\end{equation}
Define $F\hbox{:}\ \R \times \Ht \cap \Hz \longrightarrow \Lt$ by
$F(t, u) = -{\rm div}(A_t \n u)+ f \circ \phi_t \gamma_t$. Then
$F$ is a $\C^1$-map. Further $D_2F \arrowvert_{(0,y)}(0, u) =
-{\rm div}(\n u)$ (recall $y=y(\Om)$). By the standard theory of
Dirichlet boundary value problem on compact Riemannian manifolds
(\cite{1}, Theorem~4.8, p.~105 and \cite{2}, Theorem~7.32,
p.~259),
\begin{equation*}
D_2F \arrowvert_{(0,y)}\hbox{:}\ \Ht \cap \Hz \longrightarrow \Lt
\end{equation*}
is an isomorphism. By (2.2), $F(t, y^t) = 0~ \forall t$.
Proposition~2.1 now follows by the implicit function theorem.
\hfill $\Box$
\end{proof}

\vspace{.5pc}
\begin{definit}$\left.\right.$\vspace{.5pc}

\noindent $\dot{y}(\Om, V):= \left.\left(\frac{\rm d}{{\rm d}t}
y^t\right)\right\vert_{t=0} \in \Hz$ is called the (strong) {\it
material derivative of} $y$ in the direction of $V$.\vspace{.7pc}
\end{definit}

Consider $\Om^\prime \subset \subset \Om$.

\begin{propo}$\left.\right.$\vspace{.5pc}

\noindent The map $t \longmapsto y_t \arrowvert_{\Om^\prime}$ is a
$\C^1$-curve in $H^1(\Om^\prime)$ from a neighbourhood of $0$ in
$\R$ and ${\rm d}/{{\rm d}t} \arrowvert_{t=0} \left( y_t
\arrowvert_{\Om^\prime}\right) = \{ \dot{y}(\Om, V)- g(\n y,
V)\}\arrowvert_{\Om^\prime}$.\vspace{.5pc}
\end{propo}

\begin{proof}
There exists $\delta > 0 $ such that $\Om^\prime \subset
\Phi_t(\Om) \; \forall \,|t| < \delta$. Then $y_t
\arrowvert_{\Om^\prime} = y^t \circ
\Phi_{-t}\arrowvert_{\Om^\prime}\; \forall \,|t| < \delta$.
Proposition~2.2 now follows from Proposition~2.1 and
Proposition~2.38, p.~71 of \cite{5}.\hfill $\Box$
\end{proof}

\vspace{.4pc}
\begin{definit}$\left.\right.$\vspace{.5pc}

\noindent $y^\prime(\Om, V) := \dot{y}(\Om, V) - g(\n y, V) \in
\Ho$ is called the {\it shape derivative of} $y$ in the direction of
$V$.\vspace{.7pc}
\end{definit}

Consider the domain functional $J(\Om_t)$ defined by $J(\Om_t) :=
\int_{\Om_t} y_t\;{\rm d}V$ ($t \in \R$).

\vspace{.7pc}
\begin{definit}$\left.\right.$\vspace{.5pc}

\noindent The {\it Eulerian derivative} ${\rm d}J(\Om, V)$ of
$J(\Om_t)$ at $t=0$ is defined as
\begin{equation*}
{\rm d} J(\Om, V) := \lim_{t\longrightarrow 0}
\frac{J(\Om_t)-J(\Om)}{t}.
\end{equation*}
\end{definit}\vspace{.1pc}

\begin{propo}$\left.\right.$\vspace{.5pc}

\noindent The function $J(\Om_t)$ is differentiable at $t=0$ and ${\rm
d}J(\Om, V) = \int_\Om y^\prime \; {\rm d}V$.\vspace{.7pc}
\end{propo}

\begin{proof}
Let $L_V \omega$ denote the Lie derivative of $\omega$ with respect to
$V$, and $i_V \omega$ denote the interior multiplication of $\omega$
with respect to $V$. Then
\begin{equation*}
\frac{\rm d}{{\rm d}t}(\Phi_t^* \omega) \arrowvert_{t=0} =: L_V \omega =
(d\,i_V + i_V \,d)\,\omega = {\rm d}(i_V \omega)= {\rm div}(V)\,\omega.
\end{equation*}
Hence, by Propositions~2.1 and 2.2 we get
\begin{align*}
{\rm d}J(\Om, V) &= \lim_{t\longrightarrow 0} \int_\Om \left\{\frac{y^t
\,\Phi_t^* \omega -y\omega}{t}\right\} = \int_\Om
\left(\frac{\rm d}{{\rm d}t}\{y^t\,\Phi_t^* \omega\}
\right)_{\arrowvert_{t=0}}\\[.2pc]
&= \int_\Om \{\dot{y}+ y\, {\rm div}(V)\}\; {\rm d}V = \int_\Om
\{y^\prime\! + g(\n y, V)\!+y \,{\rm div}(V)\}\; {\rm d}V \\[.2pc]
&= \int_\Om y^\prime \; {\rm d}V + \int_\Om  {\rm d}(y\,i_V\omega)
= \int_\Om y^\prime\; {\rm d}V.
\end{align*}

$\left.\right.$\vspace{-2.5pc}

\hfill $\Box$\vspace{1.1pc}
\end{proof}

\begin{propo}$\left.\right.$\vspace{.5pc}

\noindent The shape derivative $y^\prime = y^\prime(\Om, V)$ is the weak
solution of the Dirichlet boundary value problem
\begin{equation}
\left. \begin{array}{r@{\ \,}l@{\ \,}l}
\La v &= &0 \quad \hbox{ on } ~\Om,\\[.2pc]
v\arrowvert_{\partial\Om} &= &-\frac{\partial y}{\partial n}\, g(V, n)
\end{array} \right\}
\end{equation}
in the space $\Ho$. {\rm (}Here{\rm ,} $n$ is the outward unit normal
field on $\partial \Om${\rm )}.
\end{propo}

\begin{proof}
Consider $\psi \in \smc$ having support in a domain $\Om^\prime \subset
\subset \Om$. There exists $\delta > 0$ such that $\Om^\prime \subset
\Om_t\; \forall \; |t|<\delta$. By problem (2.1),
\begin{equation*}
\int_{\Om^\prime} g(\n y_t, \n \psi)\;{\rm d}V = -\int_{\Om^\prime}
f\,\psi \; {\rm d}V \quad \hbox{for} \  |t|< \delta.\tag{4}
\end{equation*}
By Proposition 2.2, differentiation of LHS of eq.~(4) with respect to
$t$ at $t=0$ can be carried out under the integral sign. So we get
\begin{equation*}
\int_{\Om^\prime} g(\n y^\prime, \n \psi)\; {\rm d}V = 0.
\end{equation*}
Thus $y^\prime$ satisfies $\La y^\prime = 0 $ weakly on $\Om$.

Now $\dot{y}$, $y \in \Ht \cap \Hz$, and $y^\prime = \dot{y}-g(\n y, V)
\in \Ho$. So by Proposition~2.39, p.~88 of \cite{2}, we get
\begin{equation*}
y^\prime \arrowvert_{\partial \Om} = \dot{y}\arrowvert_{\partial
\Om}-g(\n y, V)\arrowvert_{\partial \Om} \quad \hbox{and} \quad
\dot{y}\arrowvert_{\partial \Om} =0.
\end{equation*}
Also, $y \in \smb$ and $y=0$ on $\partial \Om$ by (2.1). So, $g(\n y,
V)\arrowvert_{\partial \Om} = \frac{\partial y}{\partial n}\, g(V, n)$.
Thus, $y^\prime \arrowvert_{\partial \Om} = - \frac{\partial y}{\partial
n}\, g(V, n) $.\hfill $\Box$
\end{proof}

\section{Shape calculus for the eigenvalue problem}

\setcounter{equation}{0}
\setcounter{theore}{0}

Let $(M, g)$, $V$, $\Phi_t$, $\Om$, $\Om_t$, $\gamma_t$, $A_t$ be as in
\S2. Consider problem (1.2) posed in $\Om_t$:
\begin{equation}
\left. \begin{array}{r@{\ }l@{\quad}ll}
-\La u &= \lambda u &\hbox{on} &\Om_t , \\[.2pc]
u &= 0 &\hbox{on} &\partial \Om_t.
\end{array} \right\}
\end{equation}
Let $\lambda_1(t) := \lambda_1(\Om_t)$ and $y_1(t) := y_1(\Om_t)$ be as
in \S1. We denote $y_1(\Om)$ by $y_1$ and $\lm$ by $\lambda_1$
throughout this section.

Denote $y_1(t)\circ \Phi_t\arrowvert_\Om$ by $y_1^t$ ($t\in \R$).\vspace{.7pc}

\begin{propo}$\left.\right.$\vspace{.5pc}

\noindent The map $t \longmapsto \left( \,\lambda_1(t)\,,\,
y_1^t\,\right)$ is a $\C^1$-curve in $\R \times \Ht \cap
\Hz$ from a neighbourhood of $\,0$ {\it in} $\R$.\vspace{.7pc}
\end{propo}

\begin{proof}
By problem (3.1), for each $t \in \R$, $y_1(t)$ satisfies the equation
\begin{equation*}
\int_{\Om_t} g(\n y_1(t),\n \psi)\; {\rm d}V = \int_{\Om_t}
\lambda_1(t)\, y_1(t)\, \psi\; {\rm d}V \quad \forall \,\psi \in
H_0^1(\Om_t).\tag{5}
\end{equation*}
As in the proof of Proposition~2.1, eq.~(5) can be re-written as
\begin{equation*}
-\int_\Om {\rm div}(A_t \n y_1^t)\, \psi \; {\rm d}V = \int_\Om \lambda_1(t)\,
y_1^t \,\gamma_t \,\psi\; {\rm d}V \quad \forall \,\psi \in
H_0^1(\Om).\tag{6}
\end{equation*}
Therefore, $\,t \longmapsto \left(\, \lambda_1(t)\,,\, y_1^t\,\right)$
satisfies
\begin{equation}
\left. \begin{array}{r@{\ }lll}
{\rm div}(A_t \n y_1^t) + \lambda_1(t)\, y_1^t \,\gamma_t &= 0 \quad
\hbox{on } ~\Om, \\[.2pc]
\int_\Om \left(y_1^t\right)^2 \gamma_t \;{\rm d}V &=  1.
\end{array} \right\}
\end{equation}
Let $X := \R \times \Ht \cap \Hz$. Define $F\hbox{:}\ \R \times X
\longrightarrow \Lt \times \R$ by $F(t,\, \mu, \,u) = \left( {\rm
div}(A_t \n u) + \mu u \gamma_t, \int_\Om u^2 \,\gamma_t\;{\rm d}V -
1\right)$. Then $F$ is a  $\,\C^1$-map. Further
$D_2F \arrowvert_{(0,\, \lambda_1, \,y_1)}$ $(0,\, \mu,\, u) =
\left(\La u + \lambda_1 u + \mu y_1\;,\; 2\int_\Om y_1 u\;{\rm
d}V\right)$.\vspace{.3pc}

\begin{claim}
$D_2F \arrowvert_{(0,\, \lambda_1, \,y_1)}\hbox{:}\ \R \times \Ht \cap
\Hz \longrightarrow \Lt \times \R$ is an isomorphism.

Let $(v,\,b)\in \Lt \times \R$ be arbitrary. Consider the following
problem:
\begin{equation}
\left. \begin{array}{r@{\ }lll}
\La u + \lambda_1 u + \mu y_1 &= \,v \quad \hbox{ on } ~\Om,\\[.3pc]
2\int_\Om y_1 u \;{\rm d}V &=  \,b.
\end{array} \right\}
\end{equation}
Now by Fredholm alternative, $\La u + \lambda_1\, u = v - \mu\, y_1$ has
a solution in $\Ht \cap \Hz$ if and only if $v - \mu y_1 \perp y_1$ in
$\Lt$. So, for $\mu_0:= \int_\Om v y_1\;{\rm d}V$ there exists $u_1\in
\Ht \cap \Hz$ such that $\La u_1 + \lambda_1 u_1 + \mu_0 y_1 = v$.
Moreover, the solutions of $\La u + \lambda_1 u + \mu_0 y_1 = v$ are
of the form $u= u_1 + a y_1$, $a \in \R$. Given $b\in \R$ there exists
a unique $a_0:= b/2 -\int_\Om y_1 u_1 \;{\rm d}V \in \R$ such that
$2\int_\Om y_1 u \;{\rm d}V = b$. Put $u_0= u_1 + a_0 y_1$. Thus for
$(v,\,b)\in \Lt \times \R$ there exists a unique $(\mu_0,\, u_0)\in \R
\times \Ht \cap \Hz$ such that $D_2F \arrowvert_{(0,\, \lambda_1,
\,y_1)}(0,\, \mu_0,\, u_0) =(v,\,b)$. This proves the claim.

By (3.2), $F(t,\, \lambda_1(t),\, y_1^t) = 0 \; \forall t$.
Proposition~3.1 now follows by the implicit function theorem.\hfill
$\Box$\vspace{-.5pc}
\end{claim}
\end{proof}

\vspace{.4pc}
\begin{definit}$\left.\right.$\vspace{.5pc}

\noindent $\dot{y}_1(\Om, V) := ( ({\rm d}/{{\rm d}t}) y_1^t)
\arrowvert_{t=0} \in \Hz$ is called the (strong) {\it material
derivative of} $y_1$ in the direction of $V$.\vspace{.7pc}
\end{definit}

Consider $\Om^\prime \subset \subset \Om$.\vspace{.4pc}

\begin{propo}$\left.\right.$\vspace{.5pc}

\noindent The map $t \longmapsto y_1(t)\arrowvert_{\Om^\prime}$ is a
$\C^1$-curve in $H^1(\Om^\prime)$ from a neighbourhood of $0$ in $\R$
and $( ({\rm d}/{{\rm d}t}) [y_1(t) \arrowvert_{\Om^\prime}])
\arrowvert_{t= 0} = \left(\dot{y}_1- g(\n y_1,
V)\right)\arrowvert_{\Om^\prime}\in H^1(\Om^\prime)$. Further{\rm ,}
$y_1^\prime$ satisfies $y_1^\prime = \dot{y}_1- g(\n y_1, V)$ in $\Ho$
and $y_1^\prime \arrowvert_{\partial \Om} = - \frac{\partial
y_1}{\partial n}\; g(V, n)$.\vspace{.5pc}
\end{propo}

\begin{proof}
There exists $\delta > 0 $ such that $\Om^\prime \subset \Phi_t(\Om) \;
\forall \,|t| < \delta$. The first part of Proposition 3.2 follows from
Proposition~3.1 and Proposition~2.38, p.~71 of \cite{5}. Now as $\dot{y}_1 \in \Ho$
and $\n y_1 \in \smb$, we get $y_1^\prime = \dot{y}_1- g(\n y_1, V)\in
\Ho$. Hence, $y_1^\prime\arrowvert_{\partial \Om}=
\dot{y}_1\arrowvert_{\partial \Om}- g(\n y_1, V)\arrowvert_{\partial
\Om}=- \frac{\partial y_1}{\partial n}\; g(V, n)$.\hfill $\Box$
\end{proof}
\vspace{.4pc}
\begin{definit}$\left.\right.$\vspace{.5pc}

\noindent The {\it shape derivative of} $y_1$ in the direction of
$V$ is the element $y_1^\prime=y_1^\prime(\Om, V)\in \Ho$ defined by
$y_1^\prime =\dot{y}_1-g( \n y_1, V)$.\vspace{.7pc}
\end{definit}

\begin{propo}$\left.\right.$\vspace{.5pc}

\noindent The shape derivative $y_1^\prime\in \Ho$ satisfies
\begin{equation*}
-\La y_1^\prime = \lambda_1 y_1^\prime + \lambda_1^\prime y_1 \quad
\hbox{on } ~\Om
\end{equation*}
in the sense of distributions.
\end{propo}

\begin{proof}
Let $\psi \in \C_0^\infty(\Om)$. Let $\Om^\prime \subset \subset \Om$ be
a domain such that supp$\,\psi \subset \Om^\prime$. As $y_1(t)$ is a
solution of problem (1.2) posed in $\Om_t$, for $t$ sufficiently small we
get
\begin{equation*}
\int_{\Om^\prime} g(\n y_1(t), \n \psi)\; {\rm d}V = \int_{\Om^\prime}
\lambda_1(t)\,y_1(t)\, \psi \;{\rm d}V.\tag{7}
\end{equation*}
By Propositions 3.1 and 3.2, we can differentiate with respect to $t$
under the integral sign in eq.~(7). Thus we have
\begin{equation*}
\int_{\Om^\prime} g(\n y_1^\prime, \n \psi)\; {\rm d}V \;
= \int_{\Om^\prime} (\lambda_1 y_1^\prime + \lambda_1^\prime y_1)\, \psi
\;{\rm d}V.
\end{equation*}
Hence,
\begin{equation*}
-\int_{\Om} y_1^\prime\, \La \psi\; {\rm d}V = \int_{\Om} (\lambda_1
y_1^\prime + \lambda_1^\prime y_1)\, \psi \;{\rm d}V \quad \forall \; \psi
\in \C_0^\infty(\Om).
\end{equation*}

$\left.\right.$\vspace{-2.8pc}

\hfill $\Box$\vspace{1.3pc}
\end{proof}

\begin{propo}
\begin{equation*}
y_1^\prime \in \C^\infty(\bar{\Om}).
\end{equation*}
\end{propo}

\begin{proof}
By Proposition~3.2, $y_1^\prime = \dot{y}_1- g(\n y_1, V)$ on $\Om$.
Hence it is enough to prove that $\dot{y}_1\in \C^\infty(\bar{\Om})$.
Consider $L:= \La + \lambda_1$, a linear elliptic operator of order 2.
Then $\dot{y}_1\in\Hz$ satisfies $L\left(\dot{y}_1\right) =
L (y_1^\prime+g(\n y_1, V)) = -\lambda_1^\prime\, y_1 +
L (g(\n y_1, V))$, by Proposition~3.3. From Proposition~3.58,
p.~87 of \cite{1}, it follows that $\dot{y}_1\in \C^\infty(\bar{\Om})$.
\hfill $\Box$
\end{proof}

\begin{propo}
\begin{equation*}
\lambda_1^\prime = - \int_{\partial \Om} \left( \frac{\partial
y_1}{\partial n} \right)^2\, g (V, n)\; {\rm d}S.
\end{equation*}
\end{propo}

\begin{proof}
We write $\lambda_1^\prime =\lambda_1^\prime \int_\Om y_1^2 \; {\rm
d}V$. By Proposition~3.3, $\lambda_1^\prime = \int_\Om \{-\La
y_1^\prime- \lambda_1 y_1^\prime\}\,y_1\; {\rm d}V$. Hence by problem
(1.2) and Proposition~3.4, we get
\begin{align*}
\lambda_1^\prime = \int_\Om \{-y_1\, \La y_1^\prime+ y_1^\prime\, \La
y_1\}\; {\rm d}V &= \int_{\partial \Om} \left\{ y_1^\prime\,\frac{\partial
y_1}{\partial n}-y_1\, \frac{\partial y_1^\prime}{\partial n}\right\}\; {\rm d}S\\[.3pc]
&= \int_{\partial \Om} y_1^\prime\,\frac{\partial y_1}{\partial n}\; {\rm d}S.
\end{align*}
Now the result follows by Proposition 3.2.\hfill $\Box$
\end{proof}

\section{Proofs of Theorem~1 and Theorem~2 for $\boldsymbol{S^n}$}

\setcounter{equation}{0}

\begin{profth}
We continue with the notations of \S1 such as $r_0, r_1, \F$, and
$y(\Om), J(\Om)$ for $\Om \in \F$ for $S^n$. For $|t|<\pi$, put
$p:=(0,\ldots, 0,1)$ and $q(t)= (0, \ldots, 0, \sin t, \cos t) \in S^n$.
The Laplace--Beltrami operator $\La$ of $(S^n, \langle, \rangle)$ is
invariant under isometries of $S^n$.
So we need to study the functional $J$ only on domains $\Om(q(t)):=
B(r_1)\setminus \overline{B(q(t), r_0)}, 0 \leq |t| < r_1 - r_0$, where
$B(r_1):= B\left(p, r_1\right)$.

We define $j\hbox{:}\ (r_0 - r_1, r_1 - r_0)\longrightarrow \R $ by $j(t)=
J(\Om(q(t)))$.

Fix $t_0$ such that $0 \leq t_0 < r_1 - r_0$ and put $\Om :=
\Om(q(t_0))$ and $B_0:= B(q(t_0), r_0)$. Fix $r_2$ such that $r_0 <
r_2 < r_1 - t_0$ and consider a smooth function $\rho\hbox{:}\ S^n
\longrightarrow \R$ satisfying $\rho = 1$ on $\overline{B(q(t_0), r_2)}$
and $\rho = 0$ on $\partial B(r_1)$. Let $V$ denote the vector field on
$S^n$ defined by $V(x)= \rho(x)\, (0, \ldots, 0, x_{n+1}, -x_n)\;
\forall x=(x_1, \ldots, x_{n+1})\in S^n$. Let $\{\Phi_t\}_{t \in \R}$
be the one-parameter family of diffeomorphisms of $S^n$ associated with
$V$. Then for $t$ sufficiently close to $0$, $J(\Phi_t(\Om))= j(t_0+t)$.
Note that $J(\Phi_t(\Om))=\int_{\Om_t} y_t\; {\rm d}V$, hence by
Proposition~2.3, $j$ is differentiable at $t_0$.

Note that $j$ is an even function which is differentiable at $0$. Hence
$j^\prime(0)=0$.

Now onwards we fix $t_0$ such that $0<t_0<r_1-r_0$ and consider $\Om :=
\Om(q(t_0))$ and $B_0:= B(q(t_0), r_0)$. Let $n$ denote the outward unit
normal of $\Om$ on $\partial \Om$. For $x \in \partial B_0$, put $a=
d(p, x)$ and $\alpha = $ the angle at $p$ of the spherical triangle $T:=
[p, q(t_0), x]$ with vertices $p,\, q(t_0)$ and $x$. Then $n(x)= (q(t_0)
- \cos r_0 \,x)/ \sin r_0$ and $\left<V, n\right> (x)= (\cos a \,\sin
t_0 - \sin a \,\cos t_0 \,\cos \alpha)/ \sin r_0$. Hence, by
eq.~(19) on p.~30 of \cite{6}, we get
\begin{equation*}
\left \langle V, n\right\rangle (x)= \cos \beta (x),\tag{8}
\end{equation*}
where $\beta (x)$ denotes the angle at $q(t_0)$ of the spherical
triangle $T$ defined above.

By Proposition 2.3, $j^\prime (t_0)= \int_\Om y^\prime \; {\rm d}V$.
Hence by Proposition 2.4 and problem (1.1),
\begin{align*}
\int_\Om y^\prime \, {\rm d}V=-\int_\Om \left\{y^\prime \La y - y\, \La
y^\prime \right\}{\rm d}V &= -\int_{\partial \Om} \left\{
y^\prime\,\frac{\partial y}{\partial n}-y\, \frac{\partial
y^\prime}{\partial n}\right\} {\rm d}S\\[.3pc]
&= -\int_{\partial \Om} y^\prime\,\frac{\partial y}{\partial n}\, {\rm
d}S.
\end{align*}
Again by Proposition~2.4 and eq.~(8) above, we get
\begin{equation*}
j^\prime (t_0)= \int_{x \in \partial B_0} \left(\frac{\partial
y}{\partial n}(x)\right)^2 \cos \beta(x) \; {\rm d}S. \tag{9}
\end{equation*}
Let $H$ denote the hyperplane in $\R^{n+1}$ through $(0, \ldots, 0)$
having $q^\prime(t_0)$ as a normal vector. Let $r_H$ denote the
reflection of $S^n$ about $H$. Put $\mathcal {O} = \{ x \in \Om\, | \,
\langle x, q^\prime(t_0)\rangle >0\}$. Then $r_H(\mathcal {O}) \subset
B(r_1)$ and $r_H(\overline{B_0}) =\overline{B_0}$. For $x \in \partial
B_0 \cap \partial \mathcal {O}$, let $x^\prime$ denote $r_H(x)$. Note
that for all $x \in \partial B_0 \cap \partial \mathcal {O}$, $\cos
\beta (x) < 0$ and $\cos \beta (x^\prime) = -\cos \beta (x)$. Thus
eq.~(9) can be re-written as
\begin{equation*}
j^\prime (t_0)= \int_{x \in \partial B_0\cap \partial \mathcal {O}}
\left\{\left(\frac{\partial y}{\partial
n}(x)\right)^2-\left(\frac{\partial y}{\partial
n}(x^\prime)\right)^2\right\}\; \cos \beta(x) \; {\rm d}S. \tag{10}
\end{equation*}
The Laplace--Beltrami operator $\La$ of $S^n$ is uniformly elliptic on
$S^n$ and hence the maximum principle (\cite{4}, Theorem~5, p.~61) and
the Hopf maximum principle (\cite{4}, Theorem~7, p.~65) are applicable
on $\bar{\Om}$. Hence, by arguments analogous to \cite{3} at this stage,
we get
\begin{equation*}
\left| \frac{\partial y}{\partial n}(x)\right| < \left|\frac{\partial
y}{\partial n}(x^\prime)\right|\quad \forall \;x \in \partial B_0\cap
\partial \mathcal {O}.
\end{equation*}
Thus from eq.~(10), $j^\prime (t_0)> 0$. This completes the proof of
Theorem~1 for $S^n$.\hfill $\Box$\vspace{.7pc}
\end{profth}

\vspace{.1pc}
\begin{profthe}
We continue with the notations of \S1 such as $\lm, \,y_1(\Om)$
and $J_1(\Om)$ for $\Om \in \F$. Let $p,\, q(t)$ be as in the proof of
Theorem~1. Define $j_1\hbox{:}\ (r_0 - r_1,\, r_1 - r_0)\longrightarrow
\R$ by $j_1(t)= J_1(\Om(q(t)))$. As in the proof of Theorem~1, fix $t_0$
such that $0 \leq t_0 < r_1 - r_0$ and put $\Om := \Om(q(t_0))$ and
$B_0:= B(q(t_0), r_0)$. Then for $t$ sufficiently close to $0$ we have
$j_1(t_0+t)=J_1(\Phi_t(\Om))= \lambda_1(\Phi_t(\Om))$. By
Proposition~3.1, $j_1$ is differentiable at $t=t_0$ and $j_1^\prime
(t_0) = \lambda_1^\prime(\Om)$. As $\lambda_1(\Phi_t(\Om))
=\lambda_1(\Phi_{-t}(\Om))$, $j_1$ is an even function which is
differentiable at $0$. Thus $j_1^\prime(0)=0$.

Now onwards we fix $t_0$ such that $0<t_0<r_1-r_0$ and put $\Om :=
\Om(q(t_0))$ and $B_0:= B(q(t_0), r_0)$. Then by Proposition~3.5 and
eq.~(8), we get
\begin{equation*}
j_1^\prime (t_0) = \lambda_1^\prime(\Om) =-\int_{\partial \Om}\left(\frac{\partial
y_1}{\partial n} \right)^2 \left<V, n\right>\,{\rm d}S =-\int_{\partial
B_0}\left(\frac{\partial y_1}{\partial n} \right)^2\cos \beta (x) \,
{\rm d}S.\tag{11}
\end{equation*}
As in the proof of Theorem~1, eq.~(11) can be re-written as
\begin{equation*}
j_1^\prime (t_0)= -\int_{x \in \partial B_0\cap \partial \mathcal {O}}
\left\{\left(\frac{\partial y_1}{\partial
n}(x)\right)^2-\left(\frac{\partial y_1}{\partial
n}(x^\prime)\right)^2\right\}\; \cos \beta(x) \; {\rm d}S.\tag{12}
\end{equation*}

The Laplace--Beltrami operator $\La$ of $S^n$ is uniformly elliptic on
$S^n$. So, the Hopf maximum principle (\cite{4}, Theorem~7, p.~65) and
the generalised maximum principle (\cite{4}, Theorem~10, p.~73) are
applicable on $\bar{\Om}$. Hence, by arguments analogous to \cite{3} we\break
get
\begin{equation*}
\left| \frac{\partial y_1}{\partial n}(x)\right| < \left|\frac{\partial
y_1}{\partial n}(x^\prime)\right|\quad \forall \;x \in \partial B_0\cap
\partial \mathcal {O}.
\end{equation*}
It follows from eq.~(12) that $j_1^\prime (t_0)< 0$. The proof of
Theorem~2 is now complete for\break $S^n$.\hfill $\Box$\vspace{.4pc}
\end{profthe}

\begin{rempro}
For $t\in \R$, define $q(t)= (0, \ldots, 0, \sinh t, \cosh t) \in
\Hy^n$. Put $p:=q(0)$ and $q:=q(t_0)$ ($t_0>0$). Define the vector field
$V$ on $\Hy^n$ by $V(x)= \rho(x)\, (0, \ldots, 0, x_{n+1}, x_n)\;
\forall x=(x_1, \ldots, x_{n+1})\in \Hy^n$, where $\rho\hbox{:}\ \Hy^n
\longrightarrow \R$ is as in the proof of Theorem~1 for $S^n$.

Let $n$ denote the inward unit normal of $B(q,r_0)$ on $\partial
B(q,r_0)$. Then,
\begin{equation*}
n(x) = (q - \cosh r_0 \,x)/ \sinh r_0
\end{equation*}
and
\begin{equation*}
\left<V, n\right> (x)= (x_{n+1} \,\sinh t_0 - x_n \,\cosh t_0 \,)/ \sinh
r_0 = \cos \beta(x),
\end{equation*}
where $\beta (x)$ denotes the angle at $q$ of the hyperbolic triangle
$[p, q, x]$ with vertices $p,\, q$ and $x$.

Now Theorems~1 and 2 for the hyperbolic case can be proved using shape
calculus of \S\S2 and 3 as in the case of sphere.\hfill$\Box$
\end{rempro}

\section*{Acknowledgement}

The work of the first author is supported by research scholarship from the
National Board for Higher Mathematics NBHM/RAwards.2001/643(1).

\end{document}